\documentclass[letterpaper, 10 pt, conference]{ieeeconf}  

\IEEEoverridecommandlockouts                              

\makeatletter
\let\NAT@parse\undefined
\makeatother
\usepackage{cite}
\newcommand{\bibfont}{\footnotesize}
\usepackage[nodots]{numcompress}

\usepackage{balance}

\usepackage{url}            
\usepackage{booktabs}       
\usepackage{amsfonts}       
\usepackage{nicefrac}       
\usepackage{microtype}      

\usepackage{afterpage}

\usepackage{graphicx}

\usepackage{epsfig}

\usepackage{xcolor}

\usepackage{times} 
\makeatletter
\@ifundefined{proof}{\usepackage{amsthm}}{}
\makeatother
\usepackage{amsmath}
\usepackage{amssymb}

\usepackage{subcaption}
\usepackage{caption}

\usepackage{bm}
\usepackage{xspace}

\allowdisplaybreaks

\usepackage{multirow}

\usepackage{makecell}

\newtheorem{theorem}{Theorem}

\usepackage{algorithm}
\usepackage{algorithmic}

\newcommand{\KalmanMPC}{\texttt{Kalman-MPC}\xspace}
 
\def\p(#1|#2){p\left(#1 \, | \, #2\right)}

\def\BibTeX{{\rm B\kern-.05em{\sc i\kern-.025em b}\kern-.08em
    T\kern-.1667em\lower.7ex\hbox{E}\kern-.125emX}}
\title{\LARGE \bf
Optimal Inferential Control of Convolutional Neural Networks
}

\author{Ali Vaziri and Huazhen Fang
\thanks{A. Vaziri and H. Fang are with the
Department of Mechanical Engineering,
        University of Kansas, Lawrence, KS 66045, USA.
        Email: {\tt\small \{alivaziri,fang\}@ku.edu}}%
}

\begin{document}

\maketitle
\thispagestyle{empty}
\pagestyle{empty}

\vspace{-50mm}

\begin{abstract}
Convolutional neural networks (CNNs) have achieved remarkable success in representing and simulating complex spatio-temporal dynamic systems within the burgeoning field of scientific machine learning. However, optimal control of CNNs poses a formidable challenge, because the ultra-high dimensionality and strong nonlinearity inherent in CNNs render them resistant to traditional gradient-based optimal control techniques.  
To tackle the challenge, we propose an optimal inferential control framework for CNNs that represent a complex spatio-temporal system, which sequentially infers the best control decisions based on the specified control objectives. This reformulation opens up the utilization of sequential Monte Carlo sampling, which is  efficient in searching through high-dimensional spaces for nonlinear inference. We specifically leverage ensemble Kalman smoothing, a sequential Monte Carlo algorithm, to take advantage of its   computational efficiency for nonlinear high-dimensional systems. Further, to harness graphics processing units (GPUs) to accelerate the computation, we develop a new sequential ensemble Kalman smoother based on matrix variate distributions. The smoother is capable of directly handling matrix-based inputs and outputs of CNNs without vectorization to fit with the parallelized computing architecture of GPUs. Numerical experiments show that the proposed    approach  is effective in controlling spatio-temporal systems with high-dimensional state and control spaces. All the code and data are available at \url{https://github.com/Alivaziri/Optimal-Inferential-Control-of-CNNs}.
\end{abstract}

\section{Introduction}
Spatio-temporal dynamic systems arise in various scientific and engineering disciplines, including flow dynamics,  heat transfer,  nuclear fusion,  and natural hazards. 
Control of such systems has attracted  decades-spanning interest, with a substantial body of research.  
Model-based control has been a prevailing approach. Its success springs from incorporating mathematical models into   feedback control design to effectively handle complex dynamics. Most  models have been  constructed using  partial differential equations (PDEs)~\cite{Li-Qi:Springer:2011}. These PDE-based models are expressive and structurally efficient, but   often demand great amounts of time and labor to derive and validate for real-world  systems.  Recently, machine learning has emerged   as a powerful approach for modeling spatio-temporal systems, thanks to its capacity  of extracting  accurate  representations of a system's behavior from    abundant data, while allowing for fast training and validation~\cite{Cuomo:JSC:2022,Wikle:AnnuRev:2023}. 
In this regard, convolutional neural networks (CNNs) have proven especially useful~\cite{Cuomo:JSC:2022,Wikle:AnnuRev:2023}. CNNs are a class of multi-layer deep learning models distinguished by the inclusion of  convolutional layers. These layers apply  convolution operations to capture spatial correlations and dependencies within data. With this characteristic, CNNs are  suitable for describing complex spatio-temporal dynamics. They have been used to solve PDEs   through the  physics-informed machine learning paradigm~\cite{Karniadakis:NRP:2021, Cuomo:JSC:2022}, and have also achieved remarkable  success in  building data-driven surrogates across numerous domains, e.g., weather forecasting, traffic modeling, and aerodynamics~\cite{Chattopadhyay:SciRep:2020,Zhang:AAAI:2017,Bhatnagar:CM:2019}.
\begin{figure}[t]
    \centering
    \includegraphics[width=\linewidth]{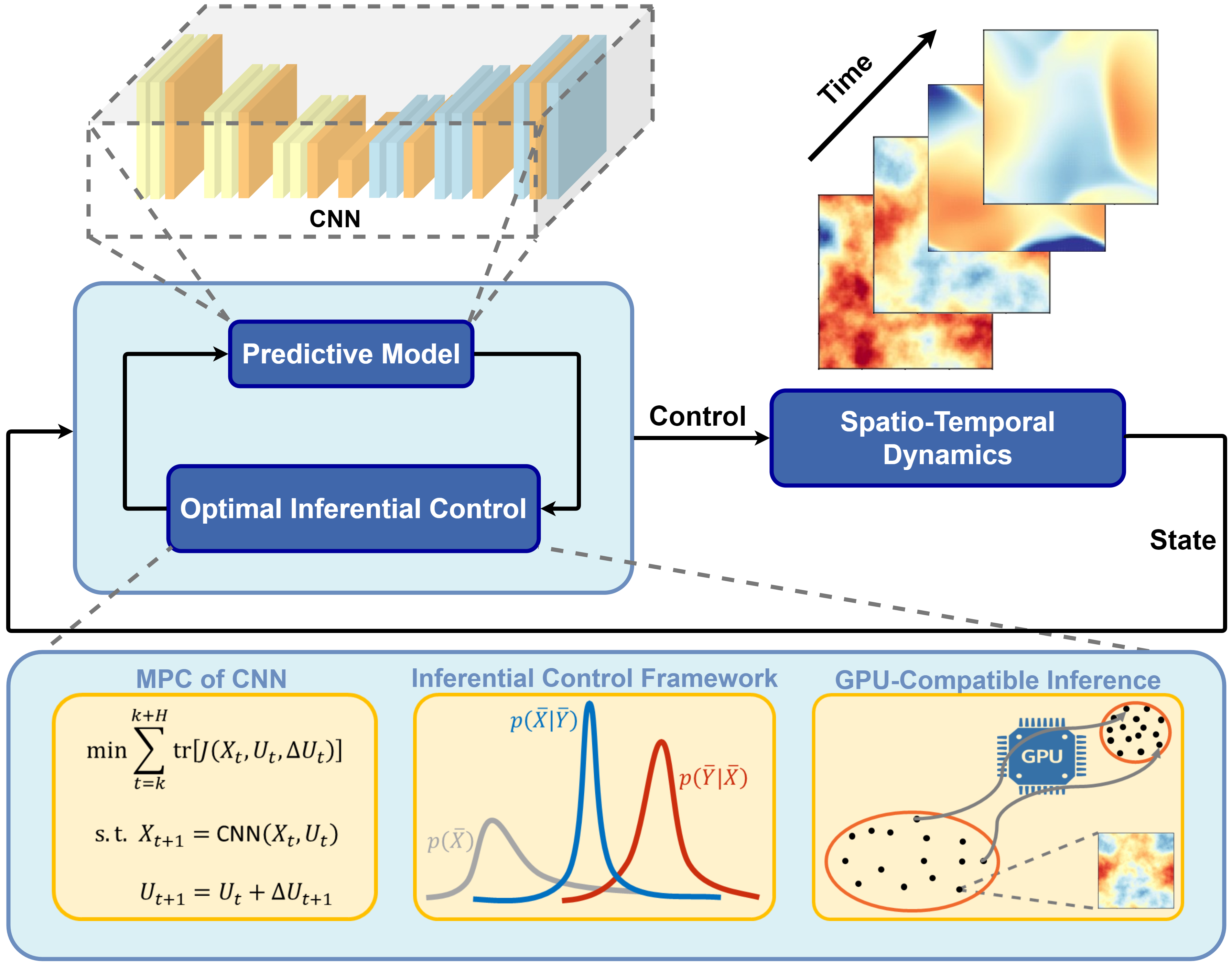}
    \caption{The OIC framework for control of CNNs. 
    }
    \label{fig:control-schematic}
\end{figure}

A CNN-based spatio-temporal dynamic model typically takes a matrix-based state-space form:
\begin{align} \label{CNN-spatio-temporal-surrogate-model}
  X_{k+1} = \mathsf{CNN}_{\theta} (  X_k,   U_k),
\end{align}
where $  X_k \in \mathbb{R}^{n \times m}$ is the spatio-temporal field at time $k$, $  U_k \in \mathbb{R}^{l \times m} $ is the control input, and $  \theta$ is the trainable parameters of the CNN. 
For~\eqref{CNN-spatio-temporal-surrogate-model}, optimal control  is an interesting new problem,   central to  applications that  will adopt CNN-based modeling. However, this problem  is  non-trivial,  because     the  ultra-high dimensions,  strong nonlinearity, and nonconvexity of $\mathsf{CNN}_{  \theta}$ will lead to extremely  complex optimization landscapes,  for which      gradient-based optimization solvers will  perform poorly or fail easily. No study has been available to investigate the challenge,  to our knowledge. 

To fill in this gap, we   develop  an {\em optimal inferential control} (OIC) framework as shown in Fig.~\ref{fig:control-schematic} to enable optimal control of CNNs, through the following two contributions. 

\begin{itemize}

\item We propose to address optimal control of CNNs via probabilistic inference. This   perspective centers around inferring or estimating  the best control actions from the control objectives and constraints, setting up the basis for the OIC framework, opening  up the use of nonlinear state estimation of~\eqref{CNN-spatio-temporal-surrogate-model}. 

\item To computationally  implement  OIC for~\eqref{CNN-spatio-temporal-surrogate-model}, we  propose  harnessing  the power of Monte Carlo sampling for inference and further combining it with the computational prowess  of graphics processing units (GPUs).  To this end, we develop a sequential ensemble Kalman smoother (EnKS). This EnKS uniquely builds on matrix variate  distributions  and   can exploit GPUs' multi-dimensional array operations to accelerate computation considerably. 

\end{itemize}

The OIC framework marks the first attempt for control of CNNs  while departing from       gradient-based   optimal control methods. It also shows the significant, yet still rarely tapped, potential of sampling-based inference and GPU to achieve computationally efficient control of highly complex systems. 

\section{Related Work}

This section provides a literature review of   the topics related with the proposed study. 

{\bf Control  of PDE systems.}   PDE control is a challenging problem that has driven many influential studies. For boundary control, the backstepping method has found important use in dealing with many types of PDEs, linear or nonlinear, and with time delays, disturbances, or parametric uncertainty~\cite{Smyshlyaev:TAC:2004, Krstic-Smyshlyaev:SIAM:2008,Coron:SIAM:2013}. Other important methods include the Riesz basis approach and Lyapunov design~\cite{GuoWang:Springer:2019,Coron:TAC:2007}. 
Optimal control of PDEs dates back to the 1970s~\cite{Lions:RMS:1973}. The works in~\cite{Bensoussan:Birkhauser:2007,Aksikas:AUTO:2009,Xu:AUTO:2011,Gahlawat:TAC:2017} investigate linear quadratic   and sum-of-squares control for   PDEs, and model predictive control (MPC) has also gained traction~\cite{Dufour:CCE:2004,Grune:SIAM:2019}. To mitigate computation, one can leverage Koopman operator theory to build data-driven linear predictors for nonlinear PDEs~\cite{arbabi_koopmanMPC, korda_koopmanMPC,  morton2018deep}, or use model order reduction to obtain low-order models~\cite{Hovland:JGCD:2008,Peitz:AUTO:2019}. 

{\bf MPC of neural network dynamic models.} The likely  first study on this topic is~\cite{Draeger:CSM:1995},   which applies dynamic matrix control, a precursor of MPC, to a neural-network-modeled pH process. Since then, growing research has led to  diverse methods and applications  spanning from chemical processes to robotics~\cite{Piche:CSM:2000,Lawryczuk:Springer:2014,Bemporad:TAC:2023,Pohlodek:ACC:2023,Cursi:RAL:2021,Salzmann:RAL:2023}. However,   gradient-based optimization, the de facto solver,   struggles to offer sufficient computation and optimization performance for deep multi-layer neural networks. To alleviate the issue,  relatively efficient solvers are custom-designed in~\cite{Salzmann:RAL:2023}, and  input-convex neural networks are used in~\cite{Bunning:L4DC:2021}. Alternatively,  sampling has shown   intriguing promise. A simple sampling technique in~\cite{Nagabandi:ICRA:2018}  is to randomly generate
many candidate control decision sequences and then pick the
sequence that leads to minimum costs. In our prior study~\cite{Askari:ACC:2021,Askari:ACC:2022}, we  propose to use particle filtering/smoothing to solve MPC of neural networks after translating it into a probabilistic inference problem; our further study~\cite{Askari:arXiv:2023} shows that improved  sampling schemes for particle filtering/smoothing can    speed up   computation greatly. Despite these advances,  CNN's ultra-high dimensionality and complicated structure present unprecedented challenges for MPC, leaving this problem widely open.

{\bf   Connections between inference and control.}
The idea of control via inference echoes some existing theories and methods. The first related concept  is the control/estimation duality, which was first established in the 1960s   between Kalman smoothing and linear quadratic regulation~\cite{Kalman:JBE:1960}. Duality for nonlinear   systems is often elusive, but the formulation of optimal control and inference problems dual to each other in this case has attracted recurring interest~\cite{Todorov:CDC:2008,Todorov:PNAS:2009,Kim:TAC:2023}. In the field of robotics, variational inference and message passing have found their way into computing optimal control problems efficiently~\cite{williams2018information,Kappen:ML:2012}.  Monte Carlo sampling, along with its variants, excels in performing probabilistic inference in   high-dimensional spaces and with complex nonlinearity. In this regard, particle filtering has shown to be useful for   MPC of nonlinear systems~\cite{Stahl:SCL:2011,Askari:ACC:2021,Askari:ACC:2022,Askari:arXiv:2023}. Reinforcement learning can also be treated as probabilistic inference~\cite{levine_arxiv_2018}, which allows to use a wide array of estimation tools to perform   learning. 

\section{Problem Setup}
 
Given~\eqref{CNN-spatio-temporal-surrogate-model}, we formulate an MPC problem in this section. First, consider the following cost function for the purpose of reference tracking: 
\begin{align*}
J\left(  X_{k:k+H},   U_{k:k+H}, \Delta  U_{k:k+H} \right) &=  \sum_{t=k}^{k+H} \mathrm{tr}  \left[   J \left(  X_t,   U_t,   \Delta  U_t \right) \right],
\end{align*} 
where $H$ is the control horizon length, and 
\begin{align*}
&  J \left(  X_t,   U_t,  \Delta   U_t \right) =   \left(   X_t -  X^\mathrm{ref}  \right)^\top   R^{-1} \left(   X_t -  X^\mathrm{ref} \right) \\
& \quad+ \left(   U_t -   U^\mathrm{ref} \right)^\top   Q_{  U}^{-1} \left(   U_t -   U^\mathrm{ref}  \right)  +  \Delta    U_{t+1}^\top    Q_{\Delta   U}^{-1}  \Delta   U_{t+1}.
\end{align*} 
Here, $  X^\mathrm{ref}$ and $  U^\mathrm{ref}$  are the references for $  X_t$ and $  U_t$, respectively, at which the model in~\eqref{CNN-spatio-temporal-surrogate-model} is at steady state satisfying  $  X^\mathrm{ref} = \mathsf{CNN}_{  \theta} \left(   X^\mathrm{ref},   U^\mathrm{ref} \right)$; $\Delta   U_t$ is the incremental control input defined as 
\begin{align} \label{incremental-input-model}
\Delta   U_{t+1} =   U_{t+1}  -   U_t.
\end{align} 
Further, $  R$, $  Q_{  U}$, and $    Q_{\Delta   U}>0$ are   weighting matrices. We now set up an MPC problem   to control the model in~\eqref{CNN-spatio-temporal-surrogate-model}: 
\begin{align}\label{original-MPC}
\begin{aligned}
   \min \quad & J\left(  X_{k:k+H},   U_{k:k+H}, \Delta   U_{k:k+H} \right) ,  \\
\mathrm{s.t.} \quad &   X_{t+1} = \mathsf{CNN}_{  \theta} (  X_t,   U_t),\\
&   U_{t+1} =   U_t + \Delta   U_{t+1},  \  t = k, \ldots, k+H. 
\end{aligned} 
\end{align} 
This problem seeks the optimal control actions within each receding horizon to steer $X_t$ towards  $X^\mathrm{ref}$ while minimizing the control effort. Its matrix-based format, which differs from conventional vector-based MPC formulations, results from controlling spatio-temporal dynamics represented by $\mathsf  {CNN}_{\theta}$. Yet, if performing vectorization, we can find  that 
\begin{align*}
&\mathrm{tr}  \left[    J \left(  X_t,   U_t,   \Delta   U_t \right) \right] = \\ 
&\quad\quad\quad   \mathrm{vec} \left(   X_t -   X^\mathrm{ref}  \right)^\top \left(   I  \otimes     R\right)^{-1}   \mathrm{vec} \left(   X_t -   X^\mathrm{ref} \right) \\
&\quad\quad\quad +  \mathrm{vec} \left(   U_t -   U^\mathrm{ref} \right)^\top \left(   I   \otimes    Q_{  U} \right)^{-1}  \mathrm{vec}\left(   U_t -   U^\mathrm{ref}  \right)\\
&\quad\quad\quad  +\mathrm{vec} \left( \Delta    U_t \right)^\top \left(   I   \otimes      Q_{\Delta   U}\right)^{-1}   \mathrm{vec} \left( \Delta   U_t\right),
\end{align*} 
which is identical to the quadratic cost  in vector-based reference-tracking MPC. Here,   $\otimes$ is the Kronecker product.  

 MPC  has been generally  addressed by gradient-based optimization. However, attempts to  apply this approach   to~\eqref{original-MPC} are   feasible only  for  small-scale CNNs. For spatio-temporal systems, a  CNN-based model often   $  X_t$ and/or $U_t$ with tens of thousands of   dimensions. Furthermore, CNNs present significant nonconvexity and nonlinearity. With these factors,  the MPC problem in~\eqref{original-MPC} is unyielding to  gradient-based optimization. As an effective solution is still absent from the literature, we propose a new framework, referred to as {\em optimal inferential control}, to tackle this challenge. Our pursuit involves transforming the problem in~\eqref{original-MPC}  into a sequential inference problem and then developing a highly efficient, GPU-compatible EnKS approach to conduct the inference.


\section{The Optimal Inferential Control Framework}

This section  derives the optimal inferential control framework for MPC of CNNs. 
To execute the framework, we develop an  EnKS approach aligned with GPU computing for computational benefits. 

\subsection{Problem Transformation}

The MPC problem in~\eqref{original-MPC} embraces a reformulation as a probabilistic inference problem. To show this, we consider a virtual auxiliary  system taking the form of
\begin{align}\label{virtual-system}
\left\{
\begin{aligned} 
  X_{t+1} &=  \mathsf{CNN}_{  \theta} \left(  X_{t},   U_{t} \right), \\
  U_{t+1} &=   U_{t} + \Delta   U_{t+1}, \\
\Delta   U_{t+1} &=    W_{t},\\
  X^\mathrm{ref} &=   X_{t}+   V_{  X, t},\\
  U^\mathrm{ref} &=   U_{t}+   V_{  U, t}, \\
\end{aligned}
\right.
\end{align} 
for $t = k,   \ldots, k+H$. In above, $  X_t$, $  U_t$, and $\Delta   U_t$ together form the virtual state,  $  X^\mathrm{ref}$ and $  U^\mathrm{ref}$ represent the virtual observation; further,  $  W_t$, $  V_{  X, t}$ and $  {V}_{  U, t}$ are additive noises following matrix normal distributions~\cite{gupta2018matrix} (see Appendix): 
\begin{gather*} 
 W_t \sim \mathcal{MN} \left(  0;   \Sigma_{  W},   \Psi_{  W} \right), \ \   V_{  X, t} \sim \mathcal{MN} \left(  0;   \Sigma_{  V_{  X}},   \Psi_{  V_{  X}} \right), \\ 
 V_{  U, t} \sim \mathcal{MN} \left(  0;   \Sigma_{  V_{  U}},   \Psi_{  V_{  U}} \right). 
 \end{gather*}
Intuitively, the virtual system in~\eqref{virtual-system} replicates the dynamics of the original system, but is observed to   behave optimally over the upcoming horizon $[k, k+H]$.  
Based on this insight, if we  estimate $  X_t$, $  U_t$, and $\Delta   U_t$ for $t = k,  \ldots, k+H$  using the virtual observation of the optimal behavior as the evidence, then the best control actions will be identified. 
To proceed forward, we rewrite~\eqref{virtual-system} compactly as
\begin{align}\label{virtual-system-compact}
\left\{
    \begin{aligned}
    \bar  {  X}_{t+1} &= f_{\mathsf{CNN}}\left(\bar {  X}_t \right) + \bar  {  W}_{t},\\
    \bar  {  Y}_{t} &=   H \bar  {  X}_t  + \bar  {  V}_{t},\\
    \end{aligned}
\right.
\end{align}
where 
\begin{gather*} 
\bar {  X}_t = 
\begin{bmatrix}
    {  X}_{t} \cr
    {  U}_{t} \cr
  \Delta   {  U}_{t} 
\end{bmatrix},  \ \ 
\bar {  Y}_t = \begin{bmatrix}
  X^\mathrm{ref}  \cr   U^\mathrm{ref}
\end{bmatrix}, \ \
\bar {  W}_t = 
\begin{bmatrix}
      0 \cr
    {  W}_{t} \cr
   {  W}_{t}  
\end{bmatrix}, \\
\bar {  V}_t = 
\begin{bmatrix}
  V_{  X, t} \cr    {V}_{  U, t}
\end{bmatrix}, \ \ 
  H = 
\begin{bmatrix}
  I &   0 &   0 \cr
  0 &   I &   0
\end{bmatrix},
\end{gather*} 
and  $f_{\mathsf{CNN}}(\cdot)$ is  evident from the context. 
To identify the state $\bar {  X}_t$ for~\eqref{virtual-system-compact}, we  take the method of maximum a posteriori  estimation, which yields
\begin{equation}\label{MAP-for-virtual-system}
   \widehat{ \bar {  X}}_{k:k+H} = \arg\max_{\bar {  X}_{k:k+H}} \log \ \p(\bar {  X}_{k:k+H}| \bar {  Y}_{k:k+H}),
\end{equation}
where $\p(\bar {  X}_{k:k+H} | \bar {  Y}_{k:k+H})$ is the  posterior distribution of 
  $\bar {  X}_{k:k+H}$ given  $\bar {  Y}_{k:k+H}$, and  $\widehat{ \bar {  X}}_{k:k+H}$ is the   estimate. This  is   a state smoothing problem in the literature on probabilistic inference. The following theorem shows an equivalence between~\eqref{MAP-for-virtual-system} and the original MPC problem in~\eqref{original-MPC}.

\begin{theorem}\label{control-inference-equivalence}
Suppose that $  W_t$, $  V_{  X, t}$ and $  {V}_{  U, t}$ are mutually independent, and that $  \Sigma_{  V_{  X}} \otimes     \Psi_{  V_{  X}} =   I \otimes   R$, $\Sigma_{  V_{  U}} \otimes     \Psi_{  V_{  U}} =  I \otimes   Q_{  U} $, and $  \Sigma_{  W} \otimes     \Psi_{  W}  =   I \otimes   Q_{  \Delta   U}$. Then,  the problems in~\eqref{original-MPC} and~\eqref{MAP-for-virtual-system} share the same optima. 
\end{theorem}

\proof
By Bayes' rule and the Markovianity   of~\eqref{virtual-system-compact}, we have 
\begin{align*}
\p(\bar { X}_{k:k+H}| \bar { Y}_{k:k+H}) &\propto \prod_{t=k}^{k+H} \p(\bar { Y}_t | \bar {  X}_t)  
 \prod_{t=k+1}^{k+H} \p(\bar {  X}_t | \bar {  X}_{t-1}) \\
 &\quad\quad \times p\left( \bar {  X}_k \right).
\end{align*} 
The log-likelihood then is given by  
\begin{align*}
&\log  p( \bar {  X}_{k:k+H} | \bar {  Y}_{k:k+H} ) \propto 
\sum_{t=k}^{k+H} \log \p(   X^\mathrm{ref} |   X_t )   \\ &\quad\quad\quad
+ \sum_{t=k}^{k+H} \log \p(  U^\mathrm{ref} |   U_t) 
+\sum_{t=k+1}^{k+H} \log p( \Delta   U_t),
\end{align*} 
as $\p(\bar {  Y}_t | \bar {  X}_t)  = \p(   X^\mathrm{ref} |   X_t)  \p(   U^\mathrm{ref} |   U_t)$, 
$\p(\bar {  X}_t | \bar {  X}_{t-1}) = p\left( \Delta   U_t \right)$, and $\bar {  X}_k$ is known at time $k$. Given~\eqref{virtual-system}, $  X^\mathrm{ref} \; | \;   X_t \sim \mathcal{MN} \left(   X_t;     \Sigma_{  V_{  X}},   \Psi_{  V_{  X}} \right)$, $  U^\mathrm{ref} \; | \;   U_t \sim \mathcal{MN} \left(   U_t;     \Sigma_{  V_{  U}},   \Psi_{  V_{  U}} \right)$, $\Delta   U_t \sim \mathcal{MN} \left(   0;     \Sigma_{  W},   \Psi_{  W} \right)$. 
This implies 
\begin{align*}
&\log \p(\bar {  X}_{k:k+H}| \bar {  Y}_{k:k+H}) \propto - \sum_{t=k}^{k+H} \mathrm{tr} \Big[  \Sigma_{  W}^{-1} \Delta   U_t    \Psi_{  W}^{-1} \Delta   U_t ^\top \\ 
& \quad\quad\quad  +    \Sigma_{  V_{  X}}^{-1} (   X_t -   X^\mathrm{ref} )   \Psi_{  V_{  X}}^{-1}  (    X_t -   X^\mathrm{ref} )^\top  \Big.\\
& \quad\quad\quad  \Big. +   \Sigma_{  V_{  U}}^{-1} (   U_t -   U^\mathrm{ref} )   \Psi_{  V_{  U}}^{-1}  (    U_t -   U^\mathrm{ref} )^\top \Big]\\
& =  - \sum_{t=k}^{k+H} \Big[ \mathrm{vec} (   X_t -   X^\mathrm{ref}  )^\top (   \Sigma_{  V_{  X}} \otimes     \Psi_{  V_{  X}} )^{-1} \mathrm{vec} (   X_t -   X^\mathrm{ref} ) \\
& \quad\quad\quad  +  \mathrm{vec} (   U_t -   U^\mathrm{ref} )^\top (   \Sigma_{  V_{  U}} \otimes     \Psi_{  V_{  U}}  )^{-1} \mathrm{vec}(   U_t -   U^\mathrm{ref}  )  \\
& \quad\quad\quad    +  \mathrm{vec} ( \Delta    U_t )^\top (    \Sigma_{  W} \otimes     \Psi_{  W}  )^{-1} \mathrm{vec} ( \Delta   U_t) \Big],
\end{align*} 
which is the opposite of  $J\left(  X_{k:k+H},   U_{k:k+H}, \Delta   U_{k:k+H} \right)$ in the problem~\eqref{original-MPC} when $  \Sigma_{  V_{  X}} \otimes     \Psi_{  V_{  X}} =   I \otimes   R$, $  \Sigma_{  V_{  U}} \otimes     \Psi_{  V_{  U}} =  I \otimes   Q_{  U} $, and $  \Sigma_{  W} \otimes     \Psi_{  W}  =   I \otimes   Q_{  \Delta   U}$. This completes the proof. \hfill \QEDclosed

Theorem~\ref{control-inference-equivalence} illustrates the viability of treating nonlinear MPC as probabilistic state estimation. This perspective not only resonates with the well-known control/estimation duality, but also sets a basis for performing MPC  via nonlinear estimation. Next, we will translate the perspective into the development of a state estimation method to deal with MPC for CNNs. Before moving forward, a   special case of Theorem~\ref{control-inference-equivalence} is noteworthy---the equivalence between ~\eqref{original-MPC} and~\eqref{MAP-for-virtual-system}   holds  when  $    \Psi_{  V_{  X}} =    R$, $  \Psi_{  V_{  U}} =    Q_{  U} $, and $    \Psi_{  W}  =   Q_{  \Delta   U}$, and $  \Sigma_{  V_{  X}}$,  $  \Sigma_{  V_{  U}}$, and $  \Sigma_{  W}=I$. This  knowledge can help one  readily set up a state estimation problem that is consistent with the    MPC problem in~\eqref{original-MPC}. 


\subsection{MPC via GPU-Compatible EnKS}

While  various nonlinear state smoothing methods exist,    they will   face  tremendous computational difficulties  if  applied directly to the problem in~\eqref{MAP-for-virtual-system}, due to the ultra-high dimensions of the CNN model. We will address this  issue in a twofold way. First, we choose EnKS to address the smoothing problem in~\eqref{MAP-for-virtual-system}.   EnKS is arguably the most computationally efficient smoothing method for high-dimensional nonlinear systems~\cite{katzfuss2019ensemble, evensen2000ensemble}. This merit stems from its design that integrates sampling-based Monte Carlo simulation with Kalman-type update in estimation. Specifically,  a   single-pass EnKS technique proposed in~\cite{evensen2000ensemble} is considered here, which computes smoothed estimates sequentially in just a single forward pass. This compares to other methods requiring two passes, forward filtering and backward smoothing, and causing more tedious computation. Second, we intend to exploit GPUs to achieve rapid computation for smoothing. GPUs are designed to support parallel computation and especially efficient in handling matrix operations. We thus will  transform the single-pass EnKS to make it compatible with the matrix-format   model in~\eqref{virtual-system-compact} and fit GPUs. 

The idea of  single-pass smoothing to deal with the   problem in~\eqref{MAP-for-virtual-system} pursues a sequential update from $\p(\bar {  X}_{k:t-1} | \bar {  Y}_{k:t-1} )$ to  $\p(\bar {  X}_{k:t} | \bar {  Y}_{k:t} )$ for $t=k, \ldots, k+H$, governed by
\begin{align}\label{single-pass-smoothing}
\p(\bar {  X}_{k:t} | \bar {  Y}_{k:t} ) = \p( \bar {  Y}_t | \bar {  X}_t) \p( \bar {  X}_t | \bar {  X}_{t-1}) \p(\bar {  X}_{k:t-1} | \bar {  Y}_{k:t-1} ).
\end{align} 
However, there is no closed-form solution to~\eqref{single-pass-smoothing} for a nonlinear system. To make the problem tractable, we define  $  {\mathcal X}_t = \bar {  X}_{k:t} $ and $  {\mathcal Y}_t = \bar {  Y}_{k:t} $ for notational simplicity, and approximate   $\p(  {\mathcal X}_t,  \bar {  Y}_t |   {\mathcal Y}_{t-1})$  as a matrix normal distribution: 
\begin{align} \label{matrix-normal-distribution-assumption}
\begin{bmatrix}
  {\mathcal X}_{t | t-1} \\ \bar{  Y}_{t | t-1}
\end{bmatrix} 
\sim
\mathcal{MN}\left(
\begin{bmatrix}
\widehat{  {\mathcal X}}_{t | t-1} \\ 
\widehat{\bar{  Y}}_{t|t-1}
\end{bmatrix}, 
\begin{bmatrix}
  \Sigma_{t|t-1}^{  {\mathcal X}} &   \Sigma_{t|t-1}^{  {\mathcal X \bar{  Y}}} \\
\left(   \Sigma_{t|t-1}^{  {\mathcal X \bar{  Y}}} \right)^\top &   \Sigma_{t|t-1}^{\bar{   Y}}
\end{bmatrix}
\otimes   \Psi
\right),
\end{align} 
where $\widehat{  {\mathcal X}}$ and $\widehat{\bar{  Y}}$ are the  means of $  {\mathcal X}$ and $\bar{  Y}$, and $  \Sigma$ are covariance matrices. 
Given~\eqref{single-pass-smoothing}-\eqref{matrix-normal-distribution-assumption}, the marginal conditional probability  distribution $\p(  {\mathcal X}_t  |   {\mathcal Y}_{t})$ is also   matrix-normal:
\begin{align} \label{marginal-matrix-normal}
 \left.   {\mathcal X}_t  \; \right| \;     {\mathcal Y}_t
\sim
\mathcal{MN}\left( \widehat{  {\mathcal X}}_{t | t},
  \Sigma_{t|t}^{  {\mathcal X}} \otimes   \Psi \right),
\end{align}
where
\begin{subequations}\label{Kalman-update}
\begin{align} \label{Kalman-update-mean}
 \widehat{  {\mathcal X}}_{t | t} &= \widehat{  {\mathcal X}}_{t | t-1}  +    \Sigma_{t|t-1}^{  {\mathcal X \bar{  Y}}}  \left(   \Sigma_{t|t-1}^{\bar{   Y}} \right)^{-1} 
\left(
\bar {  Y}_t - \widehat{ \bar {  Y}}_{t|t-1} 
\right),\\ \label{Kalman-update-cov}
  \Sigma_{t|t}^{  {\mathcal X}} &=   \Sigma_{t|t-1}^{  {\mathcal X}} +   \Sigma_{t|t-1}^{  {\mathcal X \bar{  Y}}}  \left(   \Sigma_{t|t-1}^{\bar{   Y}} \right)^{-1} \left(   \Sigma_{t|t-1}^{  {\mathcal X \bar{  Y}}} \right)^\top.
\end{align}
\end{subequations}
This shows a  closed-form Kalman-type smoothing update  of $\p(  {\mathcal X}_t |   {\mathcal Y}_t)$ in the matrix-normal setting. By applying it recursively through time, we will get a single-pass smoother. 

Next, to implement~\eqref{Kalman-update}, we take a Monte Carlo approach, using an ensemble of samples to approximate the probability distributions of interest and then performing the update over the samples. In this line, we let  $\p(  {\mathcal X}_t |   {\mathcal Y}_{t-1})$ be approximately represented by  the samples $\left\{    {\mathcal X}_{t|t-1}^i, i=1,\ldots,N \right\}$. The sample mean  for $\p(  {\mathcal X}_t |   {\mathcal Y}_{t-1})$  is then given by 
\begin{align} \label{mathcalX-prediction-sample-mean}
\widehat {  {\mathcal X}}_{t |t-1} &= {1 \over N} \sum_{i=1}^N  {  {\mathcal X}}_{t|t-1}^i, 
\end{align} 
Further, we use $   {\mathcal X}_{t|t-1}^i$ to compute the samples that weakly approximate  $\p(\bar {  Y}_t |   {\mathcal Y}_{t-1})$:
\begin{align}\label{Ybar-samples}
\bar {  Y}_{t|t-1}^i =   H  \bar {  X}_{t|t-1}^i  + \bar {  V}_t^i,
\end{align}
for $i=1, \ldots, N$, where $\bar {  V}_t^i$ is drawn from the distributions of $  {  V}_{  X, t}$ and $  {  V}_{  U, t}$. Then, we   compute  
\begin{subequations}\label{Ybar-sample-mean-cov}
\begin{align}
\widehat {\bar {  Y}}_{t|t-1} &= \sum_{i=1}^N \bar {  Y}_{t|t-1}^i,  \\
  \Sigma_{t|t-1}^{\bar {  Y}} &= {\lambda \over N} \sum_{i=1}^N  
\left( \bar{  Y}_{t|t-1}^i  - \widehat {\bar {  Y}}_{t|t-1} \right) 
\left( \bar{  Y}_{t|t-1}^i  - \widehat {\bar {  Y}}_{t|t-1} \right) ^\top, \\
  \Sigma_{t|t-1}^{  {\mathcal X}  \bar {  Y}} &= {\lambda    \over N} \sum_{i=1}^N 
\left( {  {\mathcal X}}_{t|t-1}^i  - \widehat {  {\mathcal X}}_{t|t-1} \right)
\left( \bar{  Y}_{t|t-1}^i  - \widehat {\bar {  Y}}_{t|t-1} \right) ^\top,
\end{align}
\end{subequations} 
where $\lambda = {1 / \mathrm{tr}\left(   \Psi \right)}$. 
By~\eqref{Kalman-update-mean}, the Kalman smoothing update is   executed for ${  {\mathcal X}}_{t|t-1}^i$ for $i=1,\ldots, N$ to find out 
\begin{subequations}\label{ensemble-based-Kalman-update}
\begin{align} \label{ensemble-based-Kalman-update-mean}
\widehat{  {\mathcal X}}_{t|t} &= {1 \over N} \sum_{i=1}^N    {\mathcal X}_{t|t}^i,\\ \label{ensemble-based-Kalman-update-samplewise}
  {\mathcal X}_{t|t}^i &=    {\mathcal X}_{t|t-1}^i   +    \Sigma_{t|t-1}^{  {\mathcal X \bar{  Y}}}  \left(   \Sigma_{t|t-1}^{\bar{   Y}} \right)^{-1} 
\left(
\bar {  Y}_t - {\bar {  Y}}_{t|t-1}^i 
\right),\\ \label{ensemble-based-Kalman-update-cov}
  \Sigma_{t|t}^{  {\mathcal X}} &=  {\lambda  \over N} \sum_{i=1}^N  
\left( {  {\mathcal X}}_{t|t}^i  - \widehat {  {\mathcal X}}_{t|t} \right) 
\left( {  {\mathcal X}}_{t|t}^i  - \widehat {  {\mathcal X}}_{t|t} \right)^\top, 
\end{align}
\end{subequations} 
where  $\left\{   {\mathcal X}_{t|t}^i, i = 1, \ldots, N \right\}$  weakly approximate   $\p(  {\mathcal X}_t |   {\mathcal Y}_{t})$.  Note that $\lambda$ in $  \Sigma_{t|t-1}^{\bar {  Y}}$ and $  \Sigma_{t|t-1}^{ {\mathcal X}\bar {  Y}}$ will cancel each other in~\eqref{ensemble-based-Kalman-update-samplewise}. 
Going forward, we use  these samples to make prediction to create an ensemble for $\p(  {\mathcal{X}}_{t+1} |   {\mathcal{Y}}_t )$: 
\begin{subequations}\label{prediction-mean-cov}
\begin{align}  \label{prediction-mean}
\widehat  {\bar{  X}}_{t+1|t}    &= \sum_{i=1}^N \bar{  X}_{t+1|t}^i,\\ \label{prediction-samplewise}
\bar{  X}_{t+1|t}^i &= f_{\mathsf{CNN}} \left( \bar {  X}_{t|t}^i \right) + \bar {  W}_{t}^i, \\  \label{prediction-concatenate}
 {\mathcal X}_{t+1|t}^i &\leftarrow \left(  {\mathcal X}_{t|t}^i , \bar{  X}_{t+1|t}^i \right), \ \  
\widehat{ {\mathcal X}}_{t+1|t}   \leftarrow \left( \widehat{ {\mathcal X}}_{t|t}, \widehat{\bar{  X}}_{t+1|t}  \right), \\  \label{prediction-cov}
  \Sigma_{t+1|t}^{\bar{  X}} &=   {\lambda \over N} \sum_{i=1}^N   
( \bar{  X}_{t+1|t}^i  - \widehat { \bar{  X}}_{t+1|t} )
( \bar{  X}_{t+1|t}^i  - \widehat { \bar{  X}}_{t+1|t} )^\top, \\  \label{prediction-cross-cov}
  \Sigma_{t+1|t}^{  {\mathcal X}  \bar{  X}} &=  {\lambda \over N} \sum_{i=1}^N 
( {  {\mathcal X}}_{t+1|t}^i  - \widehat {  {\mathcal X}}_{t+1|t} ) 
( \bar{  X}_{t+1|t}^i  - \widehat { \bar{  X}}_{t+1|t} )^\top, \\  \label{prediction-cov-concatenate}
  \Sigma_{t+1|t}^{  {\mathcal X}} &= 
\begin{bmatrix}
  \Sigma_{t |t}^{  {\mathcal X}} &    \Sigma_{t+1|t}^{  {\mathcal X}  \bar{  X}}  \cr
\left(  \Sigma_{t+1|t}^{  {\mathcal X}  \bar{  X}} \right)^\top &   \Sigma_{t+1|t}^{ \bar{  X}}
\end{bmatrix},
\end{align}
\end{subequations}  
 where the sample $ \bar {  W}_{t}^i$ is   taken from the distribution of $ \bar {  W}_{t}$, and $\leftarrow$ stands for the concatenating operation.  
Note that~\eqref{prediction-mean-cov} joins with~\eqref{mathcalX-prediction-sample-mean}  to close the loop for recursion. 

The above presents  a    prediction-update EnKS procedure  to perform MPC.  The prediction is performed via~\eqref{prediction-mean-cov}, and the update  via~\eqref{Ybar-samples}-\eqref{ensemble-based-Kalman-update}. These two steps run alternately and sequentially over time to complete the smoothing in a single forward pass. We name the  algorithm as \texttt{Kalman-MPC}, summarizing it in Algorithm~\ref{alg:EnKS-MPC}. In practice, one can choose not to execute~\eqref{ensemble-based-Kalman-update-cov} and~\eqref{prediction-cov}-\eqref{prediction-cov-concatenate} to save computation and storage, if they are not interested to determine $  \Sigma_{t|t}^{  {\mathcal X}} $,  $  \Sigma_{t+1|t}^{\bar {  X}} $, and $  \Sigma_{t+1|t}^{  {\mathcal X} \bar {  X}} $, which play no role in the update.

\subsection{Discussion}

The above optimal inferential control framework  is the first solution for  optimal control of CNN-modeled spatio-temporal dynamic  systems. 
Embodying the framework is the \KalmanMPC algorithm, which attains computational tractability  for several   reasons. First,   it leverages Monte Carlo sampling to search through the ultra-high-dimensional state and control spaces, and EnKS provides a principled way to ensure the efficiency in sampling. Second,   the single-pass EnKS technique partitions a large-scale search problem into smaller sequential searches, effectively limiting the search space at every step and restraining the computation. Finally, the \KalmanMPC algorithm, by matrix-based design,  well utilizes GPU's  multi-dimensional
array operations and parallelized computing capabilities. 

\begin{algorithm}[t]\small
        \caption{\KalmanMPC: MPC of CNNs via GPU-Compatible Single-Pass EnKS}
        \label{alg:EnKS-MPC}
        \begin{algorithmic}[1]
        
        \STATE Set up the MPC problem in~\eqref{original-MPC} for the CNN-based spatio-temporal model in~\eqref{CNN-spatio-temporal-surrogate-model} 
        
         \STATE Set up the virtual system in~\eqref{virtual-system-compact}
         
        \FOR{$k=1, 2, \ldots$}
        
        \STATE Initialize  $ \bar {  X}_{k|k}^i$ for $i = 1, \ldots, N $ 
        
                \FOR{$t=k+1, \ldots, k+H$}
        
        \vspace{5pt}
        
        \item[]{ \color{gray} \textsf{// Prediction}  }
        
        \STATE Compute  $  \bar {  X}_{t|t-1}^i$ and $  \hat {\bar {  X}}_{t|t-1}$ via~\eqref{prediction-mean}-\eqref{prediction-samplewise}
        
        \STATE Obtain  $    {\mathcal X}_{t|t-1}^i $ and $ \hat{   {\mathcal X}}_{t|t-1} $ by concatenation  via~\eqref{prediction-concatenate}
        
        \STATE {\color{gray} Compute   $  \Sigma_{t|t-1}^{\bar {  X}} $,   $  \Sigma_{t|t-1}^{  {\mathcal X} \bar {  X}} $, and  $  \Sigma_{t|t-1}^{  {\mathcal X}} $ via~\eqref{prediction-cov}-\eqref{prediction-cov-concatenate}}

        \vspace{5pt}
        \item[] { \color{gray}  \textsf{//   Update}}  
        
        \STATE Compute $\bar {  Y}_{t|t-1}^i$, $\hat{\bar {  Y}}_{t|t-1}$, $  \Sigma_{t|t-1}^{\bar {  Y}}$, and $  \Sigma_{t|t-1}^{  {\mathcal X}  \bar {  Y}}$ via~\eqref{Ybar-samples}-\eqref{Ybar-sample-mean-cov}
        
        \STATE Compute $ {  {\mathcal X}}_{t|t}^i$ and $\widehat{  {\mathcal X}}_{t|t}$ via~\eqref{ensemble-based-Kalman-update-mean}-\eqref{ensemble-based-Kalman-update-samplewise}
        
        \STATE  {\color{gray} Compute $  \Sigma_{t|t}^{  {\mathcal X}}$ via~\eqref{ensemble-based-Kalman-update-cov}}
        
        \ENDFOR
        
        \STATE Extract $\widehat {\bar {  X}}_{k|k:k+H}$  from $\widehat {  {\mathcal X}}_{k+H|k+H}$, and   $\widehat {  U}_{k|k:k+H}$ from $\widehat {\bar {  X}}_{k|k:k+H}$ 
        
        \STATE Apply $\widehat {  U}_{k|k:k+H}$  to the spatio-temporal system
        
        \ENDFOR

        \end{algorithmic}

        \vspace{5pt}
        {\color{gray}Note: Lines 8 and 11 are skippable.}  
\end{algorithm} 

\begin{figure*}[t]
    \centering
    \begin{subfigure}[b]{0.8\textwidth}
        \centering        
        \includegraphics[trim={2.2cm 1.6cm 0.9cm 0.35cm}, clip,width=\textwidth]{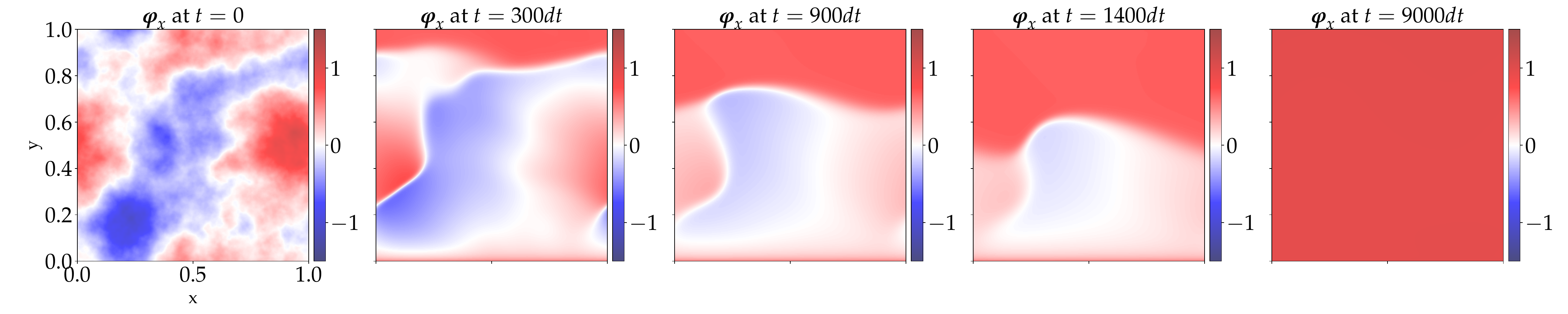}\vspace{-2mm}
        \caption{Snapshots of $  \varphi_x$ through time.}       \label{fig:phi_x_evolution_BC_Burgers}
    \end{subfigure}
    \begin{subfigure}[b]{0.8\textwidth}
        \centering        \includegraphics[trim={2.4cm 1.6cm 0.9cm 0.2cm}, clip,width=\textwidth]{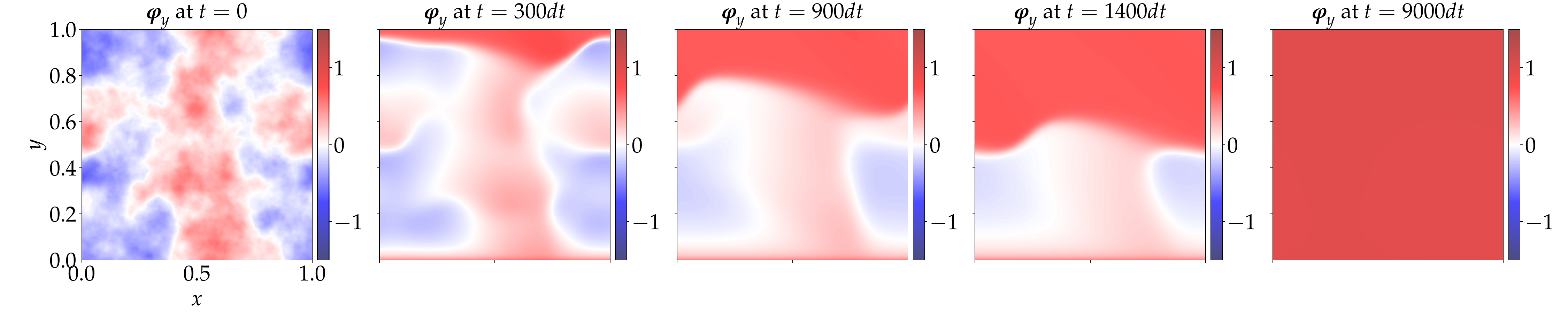}\vspace{-2mm}
        \caption{Snapshots of $  \varphi_y$ through time.}       \label{fig:phi_y_evolution_BC_Burgers}
    \end{subfigure}\\
    \begin{subfigure}[b]{0.3\textwidth}
        \centering
        \includegraphics[trim={3mm 0mm 13mm 15mm}, clip, width=\textwidth]{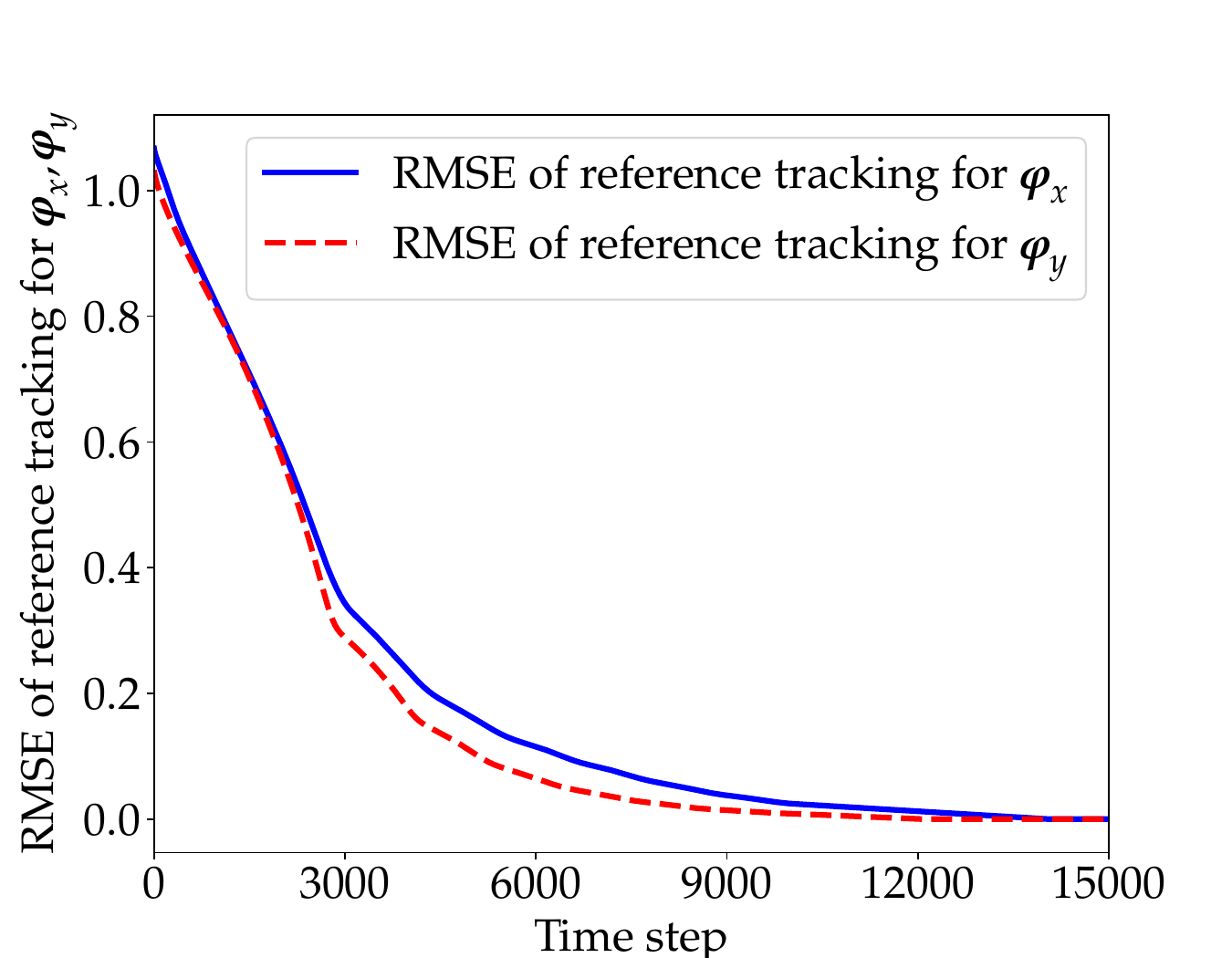}
        \caption{RMSE of reference tracking.}
        \label{fig:200by200_RMSE_uv_BC}
    \end{subfigure}
    \hspace{10mm}
    \begin{subfigure}[b]{0.3\textwidth}
        \centering
        \includegraphics[trim={3mm 0mm 13mm 15mm}, clip, width=\textwidth]{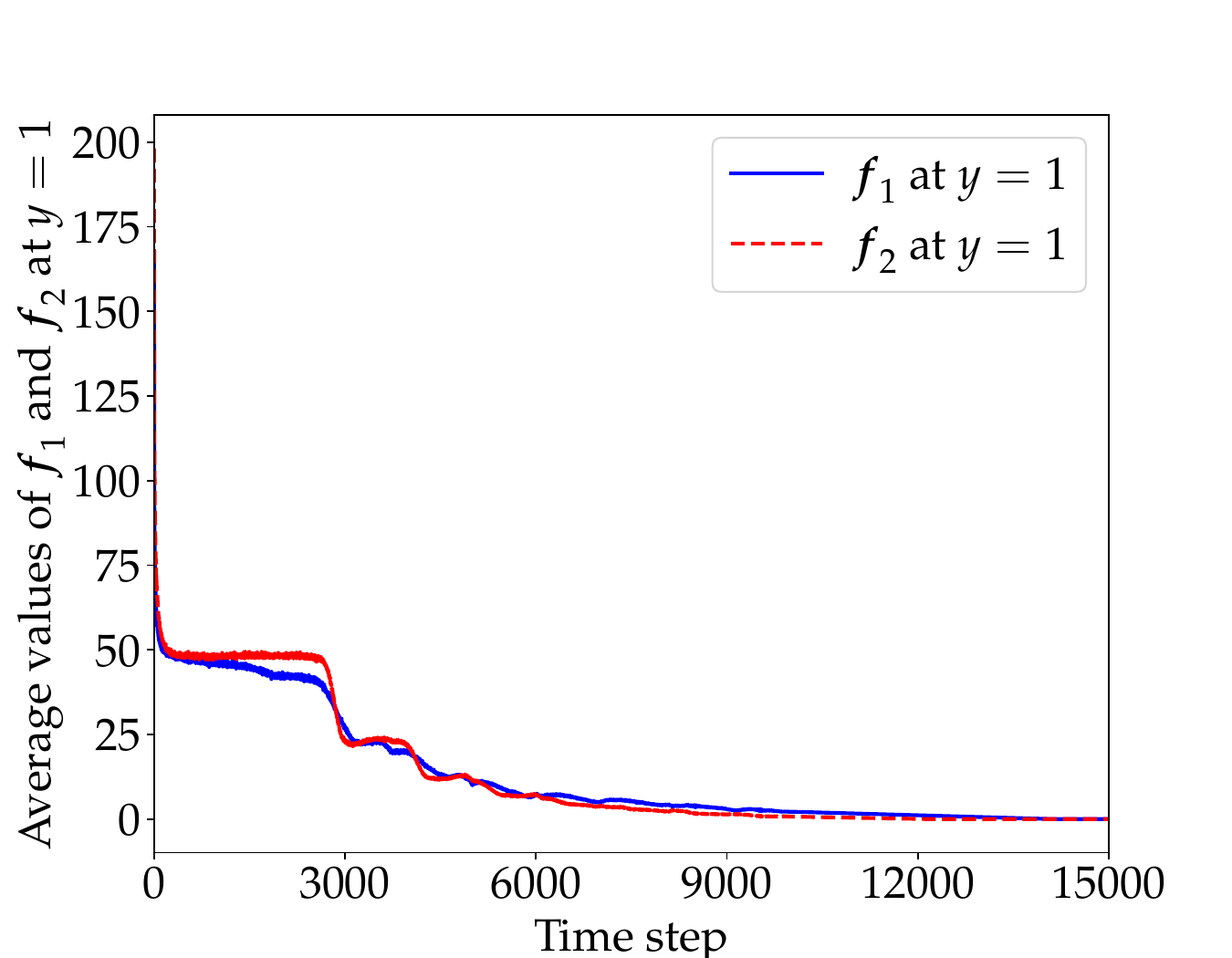}
        \caption{Control input $  f$.}
        \label{fig:plot_force1_fxfy_BC}
    \end{subfigure}   \vspace{-1mm}
    \caption{Control of the 2-D Burgers' PDE on  $200\times 200$ grid by    \KalmanMPC.} \label{fig:time_evolution_BC_Burgers}
\end{figure*} 
\begin{table*}[t]
\caption{Comparison between  the vector-based and GPU-compatible versions of   \KalmanMPC. }
\label{table_Time_Comparsion}
\small
\centering
\begin{tabular}{|c|c|c|c|c|c|} 
\hline
Grid size  &         & $2 \times 64 \times 64$ & $2 \times 80 \times 80$ & $2 \times 128 \times 128$ & $2 \times 200 \times 200$ \\ \hline
\multirow{2}{*}{\makecell{  Vector-based \\ \KalmanMPC}   }                 & Time    & $2.55$ s                & $9.51$ s                & ---                     & ---                     \\ \cline{2-6} 
 & Storage    & $26.7$ GB               & $36.1$ GB               & Memory outage             & Memory outage             \\ 
                                                 \hline
\multirow{2}{*}{\makecell{  GPU-compatible\\ \KalmanMPC}   }             & Time & $0.42$ s                & $0.51$ s                & $0.954$ s                 & $2.097$ s                 \\     \cline{2-6} 
                                                 & Storage & $2.2$ GB                & $3.2$ GB                & $4.9$ GB                  & $12.0$ GB                 \\ \hline
\end{tabular}
\end{table*}

The framework can be  expanded to include state/input constraints. Suppose supplementing the original MPC in~\eqref{original-MPC} with a constraint     $g  \left(\bar {  X}_t, \bar {  U}_t \right) \leq 0$, where $g(\cdot, \cdot)$ is   a scalar-output function  without loss of generality. We introduce a virtual observation $z_t$ to measure the constraint satisfaction:
$
z_t = \phi\left( g  \left(\bar {  X}_t, \bar {  U}_t \right)  \right) + \varepsilon_t$,  
where $ \phi(x)$ is  a barrier function taking $0$ if $x\leq 0$ and $\infty$ otherwise, and $\varepsilon_t$ is a small additive noise. By design, $z_t$ is close to $0$ if the constraint is met, and $\infty$ otherwise. We can add $z_t$ to be part of $\bar {   Y}_t$ in~\eqref{virtual-system-compact}, while setting its reference to be $0$ to ensure the awareness of the constraint in the estimation. In practical implementation, one can  use a softplus  function to make an adjustable barrier: $
 \phi(x) = \frac{1}{\alpha }\ln \left(1+e^{\beta x} \right)$,
where $\alpha, \beta > 0$. 

The  setup of $\bar {  X}_t$ requires ${  X}_t$, ${  U}_t$, and $\Delta {  U}_t$  in~\eqref{virtual-system-compact} to have the same number of columns. This is  due to the need of using matrix normal distributions and the marginalization property in~\eqref{marginal-matrix-normal}-\eqref{Kalman-update} to obtain the Kalman smoothing update in the single-pass EnKS. However,  a system's actual control input, denoted as $  u_t$, may have fewer columns. To meet the requirement, we can augment $  u_t$ to $  U_t$  by $
  U_t =   u_t   T$,
where $  T$ is a right-transformation matrix with $  T   T^\top>0$. If $  u_t$ is a column vector, a simple choice for $  T$ is $  T = \begin{bmatrix}
1 &  \cdots & 1
\end{bmatrix}$, which implies repeating $  u_t$ in the column dimension. When the \texttt{Kalman-MPC} algorithm returns $\hat  {  U}_{k|k+H}$ at time $k$, we can reconstruct  $\hat {  u}_{k|k:k+H}$ by $\hat {  u}_{k|k:k+H} = \hat  {  U}_{k|k:k+H}   T^\top  \left(   T   T^\top \right)^{-1}$.

Note that $  R$, $  Q_{  U}$, and $    Q_{\Delta   U}$  may need  some tuning  in the implementation stage, as they impact the sampling-based search within  the state and control spaces. Overall, the tuning should introduce an appropriate level of randomness to balance exploration and exploitation in the sampling. This  can be done by   adjusting   $  R$, $  Q_{  U}$, and $    Q_{\Delta   U}$  proportionately  without changing the   ratios between them. 

\begin{figure*}[t]
    \centering
    \begin{subfigure}[b]{0.8\textwidth}
        \centering        \includegraphics[trim={2cm 1.6cm 0.8cm 0.35cm}, clip,width=\textwidth]{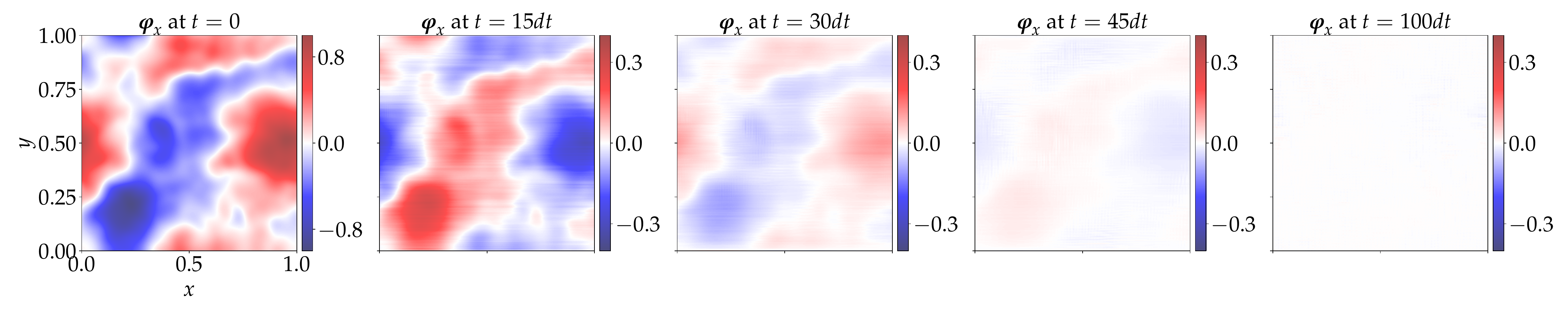}
    \end{subfigure}
    \begin{subfigure}[b]{0.8\textwidth}
        \centering        \includegraphics[trim={2cm 1.6cm 0.8cm 0.35cm}, clip,width=\textwidth]{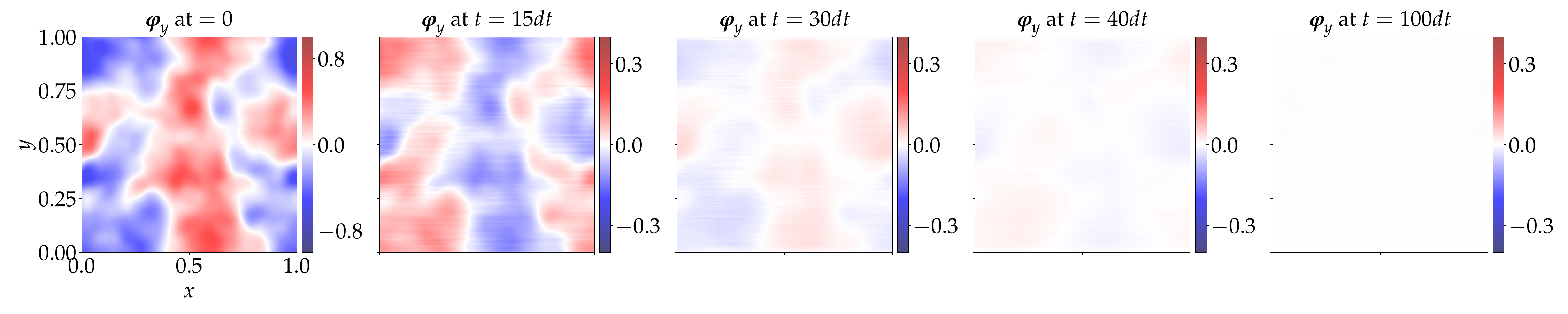}\vspace{-2mm}
        \caption{Snapshots of $  \varphi$ on $200 \times 200$ grid under control by \KalmanMPC with $100$ samples. }
        \label{phi_x-phi_y-evolution-KalmanMPC-pointwise}
    \end{subfigure}\\
    \begin{subfigure}[b]{0.8\textwidth}
    \includegraphics[trim={2cm 1.6cm 0.8cm 0.35cm}, clip,width=\textwidth]{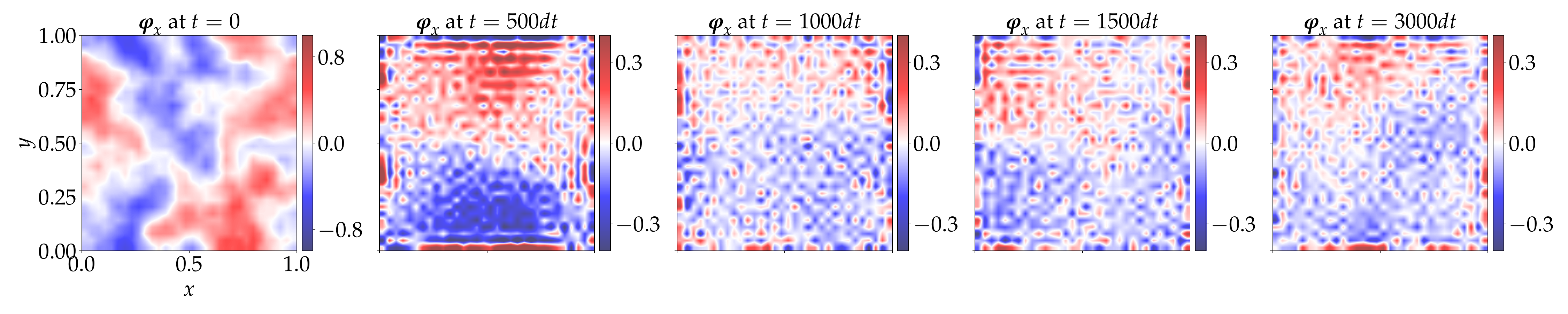}
    \end{subfigure}
    \begin{subfigure}[b]{0.8\textwidth}      \includegraphics[trim={2cm 1.6cm 0.8cm 0.35cm}, clip,width=\textwidth]{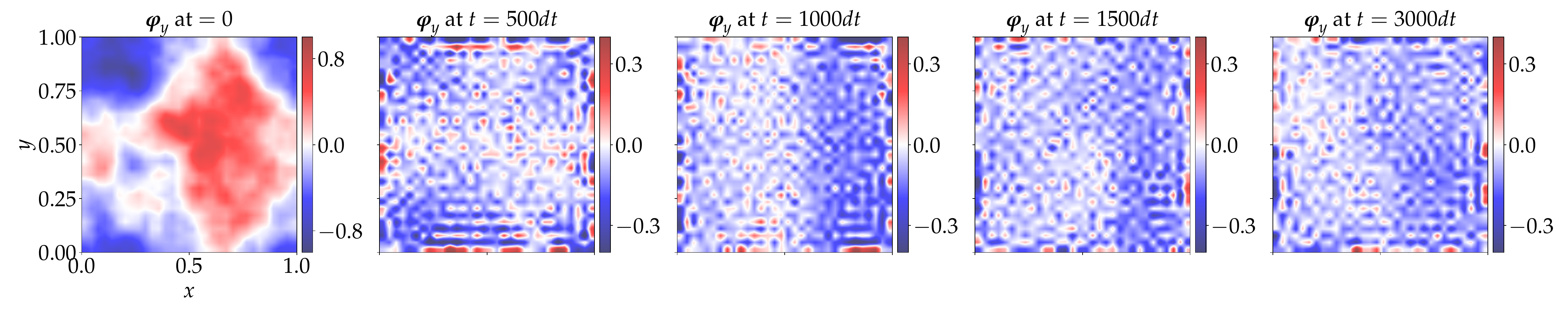}
    \vspace{-6mm}
    \caption{Snapshots of $  \varphi$  on $32 \times 32$ grid under control by IT-MPC with $1,000$ samples. }
    \label{phi_x-phi_y-evolution-IT-MPC-pointwise}
    \end{subfigure}\\
    \begin{subfigure}[b]{0.3\textwidth}
        \centering
        \includegraphics[trim={0.2cm 0 3.05cm 1.5cm}, clip, width=\textwidth]{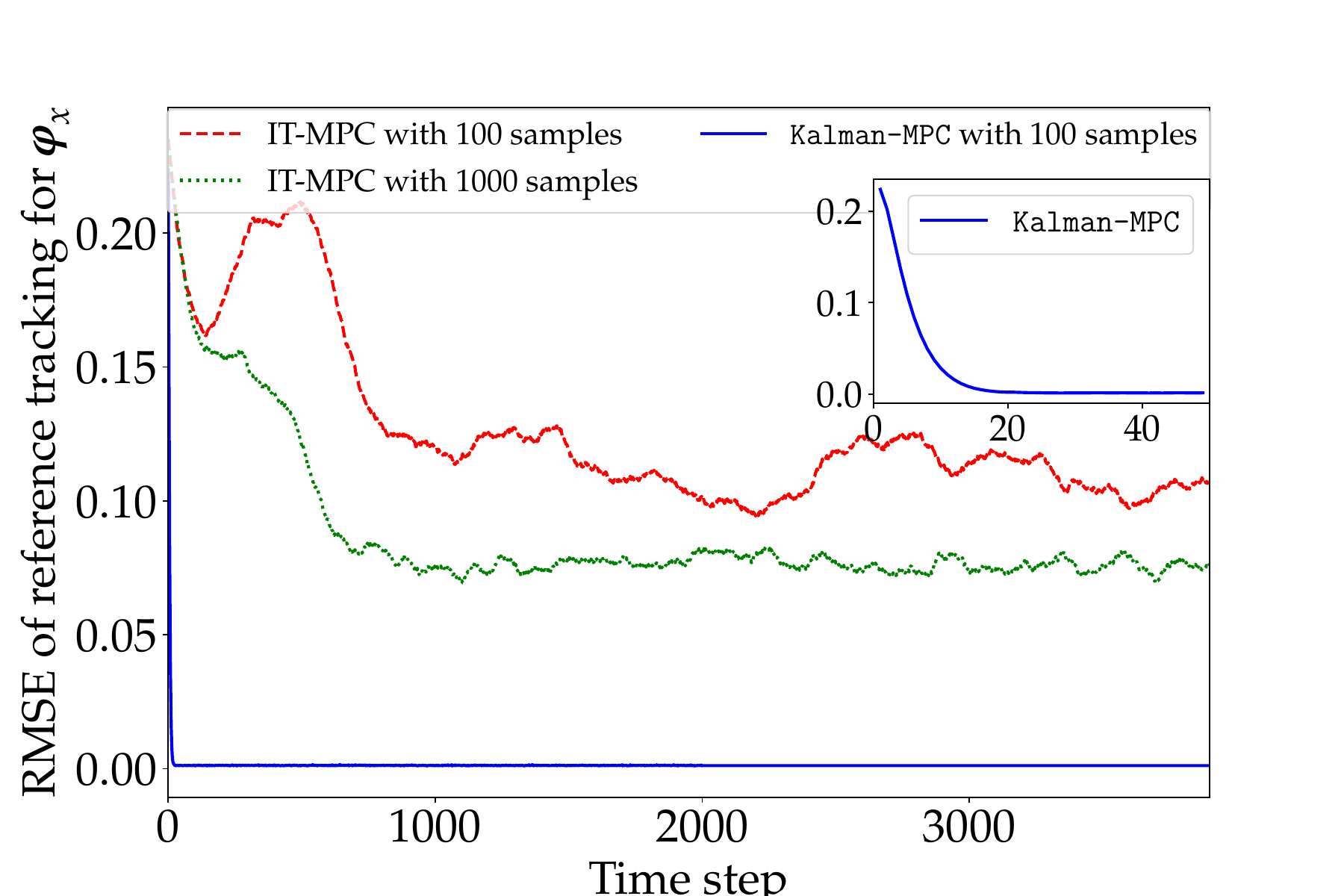}
        \caption{RMSE for $  \varphi_x$ on $32 \times 32$ grid by IT-MPC and \KalmanMPC.}
        \label{fig:MPPI_ref0_N1000_phix}
    \end{subfigure}
    \hfill
    \begin{subfigure}[b]{0.3\textwidth}
        \centering
        \includegraphics[trim={0.2cm 0 3.05cm 1.5cm}, clip, width=\textwidth]{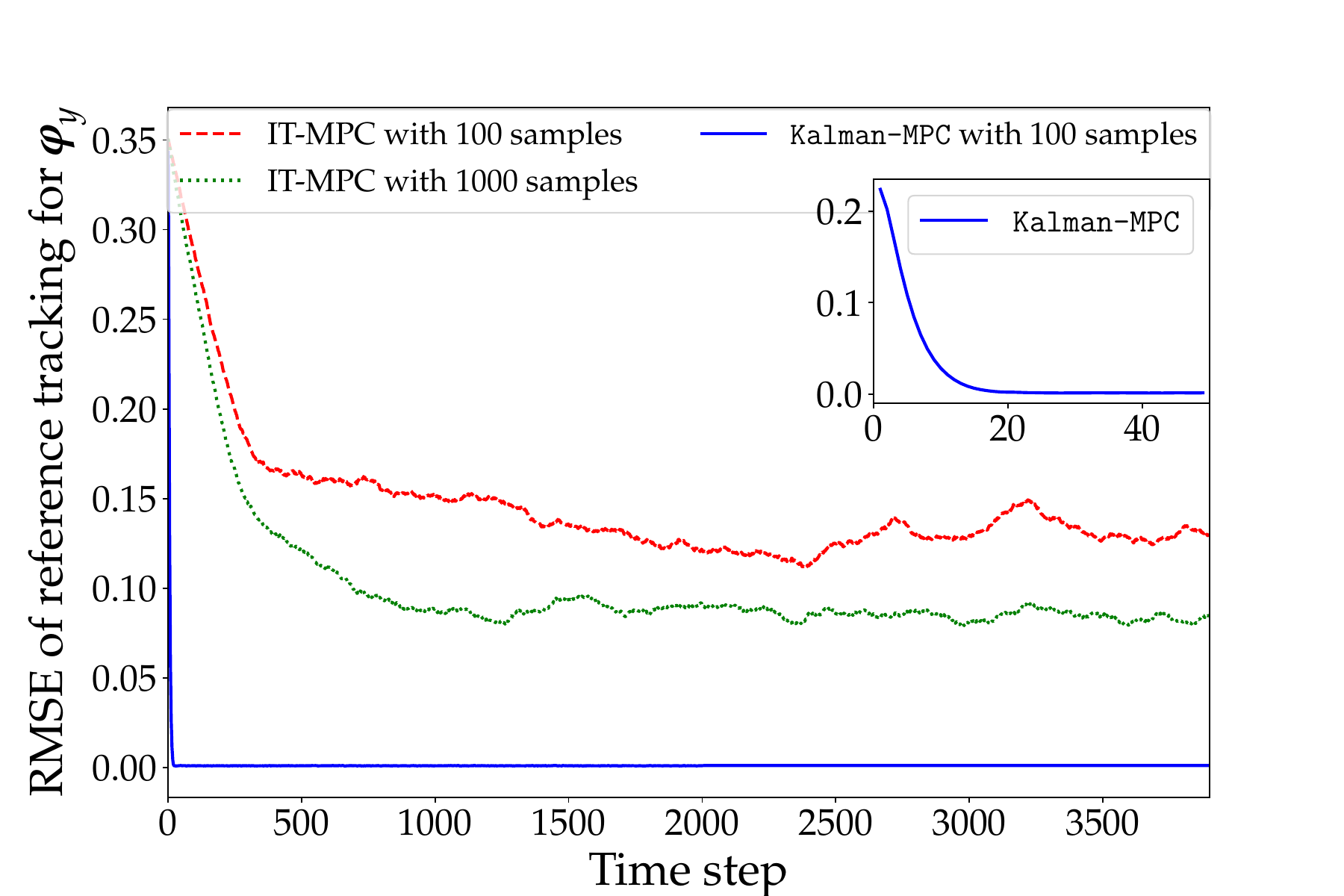}
        \caption{RMSE for $  \varphi_y$ on $32 \times 32$ grid by IT-MPC and \KalmanMPC.}
        \label{fig:MPPI_ref0_N1000_phiy}
    \end{subfigure}
    \hfill
    \begin{subfigure}[b]{0.3\textwidth}
        \centering
        \includegraphics[trim={0.2cm 0 2.3cm 1.5cm}, clip, width=\textwidth]{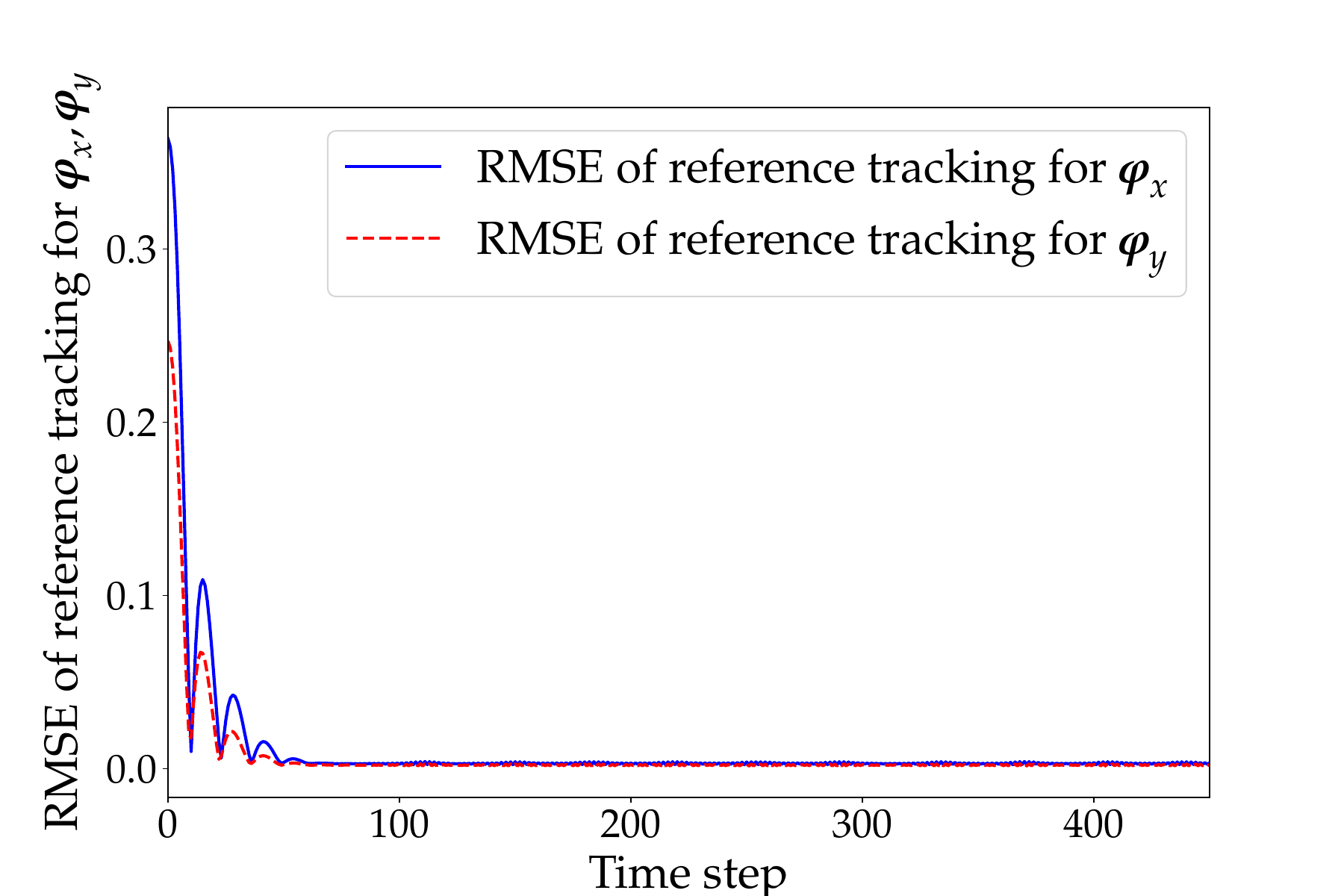}
        \caption{RMSE for $  \varphi_x$ on $32 \times 32$ grid by \KalmanMPC.}        \label{fig:200by200_pixelwise_RMSE}
    \end{subfigure}
    \vspace{-2mm}
    \caption{Comparison between \KalmanMPC and IT-MPC in control of Burgers' PDE.}    \label{fig:time_evolution_pixelwise_Burgers}
\end{figure*}

\section{Numerical Experiments}
This section evaluates   the \KalmanMPC algorithm by simulation. A CNN-based model is learned first from the 2-D Burgers' PDE, following the approach in~\cite{ren2022phycrnet}.  
We then use the  \KalmanMPC algorithm  to compute the optimal control actions based on the CNN model, and then apply the  control input to the PDE. All the simulations were conducted on  Google Colab using an NVIDIA A100-SXM4-40GB GPU and   PyTorch's CUDA library.

The 2-D Burgers' PDE is given by
\begin{align}\label{Burgers-PDE}
\frac{\partial   \varphi}{\partial t} +   \varphi \cdot \nabla   \varphi &= \nu \nabla^2   \varphi +   f, \ \ 
\end{align}
where  $  \varphi(x,y,t) = (  \varphi_x,   \varphi_y)$ is the velocity field  over the domain of $\Omega = [0, 1] \times [0,1]$, $  f$ is the external force,  $\nu=0.005$ is the viscosity coefficient. 
To train the CNN model,~\eqref{Burgers-PDE} is discretized with a $200 \times 200$ grid and solved with a time step of $dt = 0.001$. 
We consider two control setups: 1) boundary control, in which $  f$ on the boundary $y=1$, including $  f_1$ and  $  f_2$, are the control input variables, and 2) pointwise control, in which $  f$ at every spatial point is a control input variable. 
The purpose of the second setup is to test the   scalability of the \KalmanMPC algorithm to very high control dimensions.  
The \KalmanMPC algorithm runs on a horizon of $H=5$ and uses $N=100$ samples in the receding-horizon smoothing process. 
 
Beginning with the boundary control setup, we desire to steer $  \varphi$ towards the reference $  \psi^{\mathrm{ref}} =   1$, and define the root mean squared error (RMSE) in  reference tracking as \begin{align*}
\mathrm{RMSE} ={  {\left\|   \varphi -   \varphi^{\mathrm{ref}}  \right\|_F^2} \over  \mathrm{Total \ number \ of \ grid \ points}},
\end{align*}
where $\| \cdot \|_F$ is the Frobenius norm of a matrix.  Figs.~\ref{fig:phi_x_evolution_BC_Burgers}-\ref{fig:phi_y_evolution_BC_Burgers} show  the snapshots of the  velocity field at different time steps under control by the \KalmanMPC algorithm.  Initially,  $  \varphi_x$ and $  \varphi_y$ vary arbitrarily over the domain, but   both converge to  $  \psi^{\mathrm{ref}}$. 
This trend is also displayed in Fig.~\ref{fig:200by200_RMSE_uv_BC}, in which the RMSE declines to zero. The control input variables $  f_1$  and $  f_2$ also converge toward zero, as shown in Fig.~\ref{fig:plot_force1_fxfy_BC}, because $\left(   \varphi =  1,   f =   0 \right)$ is a stable equilibrium point for~\eqref{Burgers-PDE}. 

Further, we are interested to find out how much the GPU-compatible design of the \KalmanMPC algorithm improves the computation. To this end, we consider the vector-based version  of the algorithm. This variant, referred to as {\em vector-based \KalmanMPC} in sequel,   runs the regular vector-based single-pass EnKS on GPUs after vectorizing~\eqref{virtual-system-compact}.  Table~\ref{table_Time_Comparsion} offers the comparison of the two versions over different spatial grid sizes. Within the specified GPU-computation setting, the  vector-based \KalmanMPC algorithm can achieve control of Burgers' PDE  on a grid of up to $80\times 80$, but will run out of memory beyond this size. By contrast, the GPU-compatible \KalmanMPC algorithm scales well to a $200 \times 200$ grid. The different memory usage results from the different dimensions of the covariance matrices. If the grid size is $n  \times m$, $  \Sigma^{ {\mathcal X}\bar{  Y}}$ is of up to $nmH \times n$, and $  \Sigma^{\bar{  Y}}$ of $nm \times nm$  for  the  vector-based \KalmanMPC algorithm; but they are just of $nH \times n$ and $n \times n$ for    the   \KalmanMPC algorithm. Significant difference also arises in the computation time. The GPU-compatible \KalmanMPC algorithm requires much less time, by one order of magnitude, while scalable to the $200 \times 200$ grid. Our numerous experiments offline have shown that, computation-wise, this  algorithm could handle a grid of up to $500 \times 500$. But it is hard to train a sufficiently accurate CNN model for this grid, so we skip presenting the results in the paper.

Next, we consider the pointwise control setup to drive $  \varphi$ to zero. Here, the control dimensions are the same as the state dimensions,   giving more challenges for  computation. For the purpose of comparison, we consider the information-theoretic MPC (IT-MPC) algorithm proposed in~\cite{williams2018information}, which builds upon an information-theoretic interpretation of optimal control, implements  sampling-based optimization, and runs on GPUs. After many tries, we find out that the IT-MPC algorithm can only handle a grid size of at most $32 \times 32$. Fig.~\ref{phi_x-phi_y-evolution-IT-MPC-pointwise} shows the evolution of the Burger's velocity field in this case, in which the IT-MPC algorithm struggles to accomplish the control objective even though using 1,000 samples. However, the \KalmanMPC algorithm  succeeds 
on a $200\times 200$ grid using only $100$ samples, as shown in Fig.~\ref{phi_x-phi_y-evolution-KalmanMPC-pointwise}. Further, the RMSE in reference tracking as shown in Figs.~\ref{fig:MPPI_ref0_N1000_phix}-\ref{fig:200by200_pixelwise_RMSE} offers a clearer view of the difference, indicating the much better performance by the \KalmanMPC algorithm. The IT-MPC algorithm demonstrates limited  scalability, likely because it seeks to search the  state and control spaces over the entire  horizon. The \KalmanMPC algorithm, in contrast, performs sequential searches at every time step within the horizon,   requiring much fewer samples and less computation. Because of the same reason, the \KalmanMPC algorithm requires much less tuning effort to succeed than the IT-MPC does, as consistently shown in our experiments. 

\section{Conclusion}

Modeling and control of spatio-temporal dynamic systems  are  critical problems across many domains. CNNs have become increasingly popular for building data-driven models that capture complex spatio-temporal dynamics, but what has been unexplored  is    optimal control of such CNN-based  dynamic models. As a primary challenge, the ultra-high dimensions and nonlinearity intrinsic to CNNs will overwhelm numerical optimization techniques to make the computation prohibitive. In this paper, we propose the optimal inferential control framework to overcome the challenge. We convert the problem of MPC of CNNs into an equivalent probabilistic inference problem. Then, we develop a GPU-compatible EnKS to perform MPC. The resulting algorithm, \KalmanMPC, combines the potency of Monte Carlo sampling with the computing power of GPUs  to efficiently infer the best control decisions despite high-dimensional state and control spaces. Numerical experiments show that the \KalmanMPC algorithm succeeds in controlling spatio-temporal systems with tens of thousands of state and control dimensions. Our study establishes a new way to control spatio-temporal dynamic systems, with potential application to other systems involving high-dimensional nonlinear dynamics.

\section*{Appendix}

\addcontentsline{toc}{section}{Appendices}
\renewcommand{\thesubsection}{\Alph{subsection}}

\label{appendix:matrix-normal-distribution}

In this Appendix, we provide a brief overview of random matrices and the matrix normal distribution, which is taken from~\cite{gupta2018matrix}. 

A   matrix $ X \in \mathbb{R}^{n \times m}$ is a random matrix if its elements are random variables. Its probability density function (pdf) $p  \left(  X \right) $ satisfies 
\begin{itemize}
\item[1)] $p  \left(  X \right) \geq 0$, 

\item[2)] $\int_{\Omega} p  \left(  X \right) d  X = 1$, and 

\item[3)] $P\left(   X \in  S \right) = \int_{ S} p  \left(  X \right) d  X$, 
\end{itemize}
where $ \Omega$ is the space of realizations of $ X$, and $ S$ is a subset of $ \Omega$. 
If $ X \in \mathbb{R}^{n \times m}$ and $ Y \in \mathbb{R}^{r \times s}$ are two random matrices, 
their  joint pdf $p  \left(  X,  Y \right)$ satisfies
\begin{itemize}

\item[1)] $p \left(  X,  Y \right) \geq 0$, 

\item[2)] $\int_{\Omega_{ X,  Y}} p  \left(  X,   Y \right) d  X d  Y = 1$, and 

\item[3)] $P\left( (  X,    Y) \in  S \right) = \int_{ S} p  \left(  X,  Y \right) d  X  d  Y$ , 
\end{itemize}
where $ \Omega_{ X,  Y}$ is the space of realizations of $ (  X,    Y)$, and $ S$ is a subset of $  \Omega_{  X,   Y}$. The marginal pdf of $  X$ is defined by 
\begin{align*}
p \left(   X \right) = \int_{  Y} p  \left(   X,   Y \right) d   Y,
\end{align*}
and the conditional pdf is defined by
\begin{align*}
\p(   X |   Y ) = { p  \left(   X,   Y \right) \over p  \left(    Y \right)  }. 
\end{align*}

We say that $  X \in \mathbb{R}^{n \times m}$  has a matrix normal distribution if its pdf is given by
\begin{align*}
p \left(  {X} \right) &= (2 \pi)^{-\frac{nm}{2}} \left|  {\Sigma} \right|^{-\frac{n}{2}} \left|  {\Psi} \right|^{-\frac{m}{2}} \\
&\quad \times \exp \Bigg[ 
-\frac{1}{2} \mathrm{tr} \left(  {\Sigma}^{-1} \left(  {X} -  {M} \right) {\Psi}^{-1} \left(  {X} -  {M} \right)^\top \right) \Bigg],
\end{align*}
where $  M  \in \mathbb{R}^{n \times m}$, $   \Sigma \in \mathbb{R}^{n \times n}>0$, and $   \Psi \in \mathbb{R}^{m \times m}>0$. The shorthand notation is 
\begin{align*}
  X \sim \mathcal{MN} \left(   M;   \Sigma,   \Psi \right). 
\end{align*}
It is provable that $\mathrm{vec} \left(   X^\top \right)$ follows a multivariate normal distribution given by
\begin{align*}
 \mathrm{vec} \left(   X^\top \right)  \sim \mathcal{N} \left( \mathrm{vec} \left(   M^\top \right);    \Sigma \otimes   \Psi \right)  . 
\end{align*}
The following relations about expected values will hold:
\begin{align*}
\mathrm{E} \left(    X \right) &=   M,\\
\mathrm{E}\left[ \left(   X -   M \right) 
\left(   X -   M \right)^\top 
\right]&=\mathrm{tr}(  \Psi)   \Sigma,\\
\mathrm{E}\left[ \left(   X -   M \right)^\top 
\left(   X -   M \right) 
\right]&=\mathrm{tr}(  \Sigma)   \Psi. 
\end{align*}

If $  A \in \mathbb{R}^{r \times n}$ and $  B \in \mathbb{R}^{m \times s}$ are of full rank, then 
\begin{align*}
  A   X   B \sim \mathcal{MN} \left(
  A   M   B;
   A   \Sigma   A^\top,
  B^\top   \Psi   B
\right).
\end{align*}
If
\begin{align*}
\begin{bmatrix}
  X_{n \times m} \cr
  Y_{r \times m} 
\end{bmatrix}
\sim \mathcal{MN}
\left(
\begin{bmatrix}
  M_{  X} \cr
  M_{  Y}
\end{bmatrix};
\begin{bmatrix}
  \Sigma_{  X} &   \Sigma_{  X   Y}\cr
  \Sigma_{  X   Y}^\top &   \Sigma_{  Y}
\end{bmatrix},   \Psi
\right),
\end{align*}
then
\begin{align*}
 {X} \; | \;  {Y} \sim \mathcal{MN} \Big(
 {M}_{ {X}} +  {\Sigma}_{ {X}  {Y}}  {\Sigma}_{ {Y}}^{-1} \left(  {Y} -  {M}_{ {Y}} \right); \\
 {\Sigma}_{ {X}} -  {\Sigma}_{ {X}  {Y}}  {\Sigma}_{ {Y}}^{-1}  {\Sigma}_{ {X}  {Y}}^\top, 
 {\Psi}
\Big).
\end{align*}
 
\balance

\bibliographystyle{unsrt}
\bibliography{references}

\end{document}